\documentclass[12pt]{iopart}

\usepackage[compress]{cite}
\usepackage{graphicx}
\usepackage{bm}

\begin{document}

\title[Optimal networks measured by GMFRT]{Optimal networks measured by global mean first return time}

\author{Junhao Peng $^{1,2}$,  Renxiang Shao $^{1,2}$, Huoyun Wang $^{1, 2}$}
\address{1.School of Math and Information Science, Guangzhou University, Guangzhou 510006, China. \\
2.Guangdong Provincial Key Laboratory co-sponsored by province and city of Information Security Technology, Guangzhou University, Guangzhou 510006,  China.
} \ead{\mailto{pengjh@gzhu.edu.cn}}


\begin{abstract}
Random walks have wide application in real lives, ranging from target search, reaction kinetics, polymer chains, to the forecast of the arrive time of extreme events, diseases or opinions. In this paper, we consider discrete random walks on general connected networks and focus on the analysis of the global mean first return time (GMFRT), which is defined as the mean first return time averaged over all the possible starting positions (vertices), aiming at finding the   structures who have the maximal (or the minimal) GMFRT among all connected graphs with the same number of vertices and edges. Our results show that, among all trees with the same number of vertices, trees with linear structure are the  structures with the minimal GMFRT and stars are the  structures with the maximal GMFRT.  We also find that, among all connected graphs with the same number of vertices, the graphs whose vertices  have the same degree, are the  structures with the minimal  GMFRT; and the graphs whose vertex degrees have the biggest difference, are the structures with the maximal  GMFRT. We also present the methods for constructing the graphs with the maximal GMFRT (or the minimal GMFRT), among all connected graphs with the same number of vertices and edges.
\end{abstract}

05.40.Fb Random walks and Levy flights
05.60.Cd Classical transport
\maketitle

\section{Introduction}
\label{intro}

Random walks  are popular models  with wide applications~\cite{Book-Metzler-Oshanin-Redner-2014, MASUDA_PR2017}, which include  target search~\cite{Rupprecht-Benichou-2016-PRE,Cvijovic-Klinowski-1995-Science,Sims-Southall-Humphries-Hays-2009-Nature}, reaction kinetics~\cite{Ben-Avraham2000,Benichou2014-PR,Benichou-Chevalier-2010-NatChem}, descriptions of financial markets~\cite{Ghashghaie-Breymann-Peinke-Talkner-Dodge-1996-Nature, Masoliver-Montero-Perello-Weiss-2006-JEBO} and polymer chains~\cite{Cates-Witten-1987-PRA, BLUMEN200512, GurtBlumen05}. Random walks can also be used to model  epidemic or opinion  spreading\cite{NAGATANI_Physica_A2019, LI-IF-2019, Zhang_Chaos_2011, Ree2011-PRE}, help predict the arrival time of diseases ( or opinion) spreading on networks~\cite{Lannelli_KoherPRE-2017} and  estimate the occurrence (or recurrence) of extreme events on the networks~\cite{Hackl2019, BunKr05, KondVa06, BatGer02}, and etc.  Classical dynamic models show that  scale-free social networks are prone to the spreading of rumours~\cite{NEKOVEE-Physica_A-2007, GuazzCini2015}. When random-walking agents are introduced into the networks, the spreading of rumours shows different characters\cite{Zhang_Chaos_2011, Ree2011-PRE}. In conclusion, the study of random walks is  of great importance. 

One of the most attracting quantities relevant to random walks is the first passage time, which reflects how long it takes  a walker to walk from a site to the target site~\cite{Book-Redner-2007}. Results show that the structures of the networks have great effect on first passage time~\cite{PengAgliari19a, AgCasCatSar15, PengAgliariZhang15, ZhZhGa10, CampariCassi2012, PengElena2017, Wu_2019}. Another important quantity is the first return time~\cite{MaKo04,EiKaBu07,MoDa09,SaKa08,PaPen11,HadLue02,Olla07},  which is the time it takes  a random walker to return to the starting site for the first time. How fast does social opinion reach back to the sender~\cite{GuazzCini2015,TsaNi2019-FRA,Lima2017}? How long does it take two walkers, starting from the same site,  to meet again? How long is the time interval between two successive extreme events (e.g. floods, droughts, violence) in social lives~\cite{BunKr05,KondVa06,BatGer02}? All  these questions can be replied and explained by the first return time (FRT) for random walk on the corresponding  graphs. Note that the FRT is a  random variable, one can  analyze the probability distribution of  the FRT~\cite{Olla07,Ch11,MuSu10,LowMast00,IzCa06}. One can also analyze the moments of the FRT\cite{PengXu2018, liu2019scalings,PengShao_2018}. Fortunately, the  mean of the FRT (MFRT) can be exactly evaluated by using the Kac lemma~\cite{Kac1947}, and for classical discrete random walks on finite connected graphs  $G=(V, E)$ $(V=\{v_1, v_2,\cdots, v_n\})$, the MFRT for random walker starting from vertex $v_i$ $(i=1,2, \cdots, n)$ satisfies ( e.g., see~\cite{LO93, CondaminBenichou2007})
\begin{eqnarray}
\mu_i=\frac{1}{\pi_i}=\frac{2m}{d_{i}},
\label{MFRT}
\end{eqnarray}
where $\pi_i=\frac{d_{i}}{\sum_{v_j\in V}d_{j}}$, $d_{i}$ is the degree of vertex $v_i$ and $m=|E|$ is the total number of edges of graph $G$.


In order to disclose the  whole character of the MFRT for the whole graph, one can further analyze  the global MFRT (GMFRT), which is defined as
\begin{eqnarray}
\overline{\mu}&=&\frac{1}{n} \sum_{i=1}^n\mu_i=\frac{1}{n} \sum_{i=1}^n\frac{1}{\pi_i}.
      \label{Def_GMFRT}
\end{eqnarray}
  Given the same numbers of vertices and edges, which structures  have the minimal GMFRT  and which structures  have the   maximal GMFRT? The answers for these questions would be helpful for the design and optimization of networks. 
   The MFRT has close connection with the moments of first passage time\cite{PengXu2018, PengShao_2018}, and the first passage time is useful indicator for transport efficiency of networks. The results obtained in this paper would be helpful for  understanding the transport properties of the networks. 

In this paper, we study  discrete random walks on general connected networks aiming at finding the  structures  with the maximal GMFRT and the  structures with the minimal GMFRT among all connected networks with the same numbers of vertices and edges. Firstly, we consider the case while the connected networks are trees and then we analyze the GMFRT for general connected networks. For both cases, we find the structures with the maximal GMFRT (or the minimal GMFRT). Our results show that, among all trees with the same number of vertices, trees with linear structure are the  structures with the minimal GMFRT and stars are the  structures with the maximal GMFRT.  We also find that, among all connected graphs with the same number of vertices, the graphs whose vertices  have the same degree, are the  structures with the minimal  GMFRT; and the graphs whose vertex degrees have the biggest difference, are the  structures with the maximal  GMFRT. The methods for constructing the  graphs with the maximal GMFRT (or the minimal GMFRT) are also presented in this paper. 

 \section{The optimal trees among all trees  with the same number of vertices}
\label{sec:3}
Let $G=(V, E)$ be a tree with $n$ vertices. It is a connected graphs with no cycle. The total number of edges of $G$ is $m=n-1$.
Thus, for any vertex $v_i$ $(i=1,2, \cdots, n)$, the MFRT for random walk starting from  vertex $v_i$ satisfies
\begin{eqnarray}
\mu_i=\frac{1}{\pi_i}=\frac{2m}{d_{i}}=\frac{2(n-1)}{d_{i}}.
\label{EqPijC1}
\end{eqnarray}
In this section, we analyze the GMFRT of trees and find the  trees with the maximal (or the minimal) GMFRT among all trees with the same number of vertices.


\subsection{The trees with the minimal GMFRT}
\label{sec:3_1}
In this section, we find that trees with linear structure have the minimal  GMFRT among all trees with the the same number of vertices. Calculating the GMFRT, we obtain the minimal GMFRT for trees with $n$ vertices, shown as
  \begin{eqnarray}
\overline{\mu}_{min}&=&\frac{2(n-1)\times 2+ (n-1)\times (n-2)}{n} \nonumber \\
      &=& n+\frac{n-2}{n}.
\end{eqnarray}
 This result is obtained by using  the following argument: for any tree, if there is one vertex whose degree is more than $2$, we can reconstruct the tree, and obtain a new tree with less GMFRT. As a consequence, the trees, whose  vertex degrees  are less than or equal to $2$ (i.e., trees with linear structures), have the minimal  GMFRT among all trees with the same numbers of vertices. The proof of the argument is as follows.\par

\begin{figure}
\begin{center}
\includegraphics[scale=1.2]{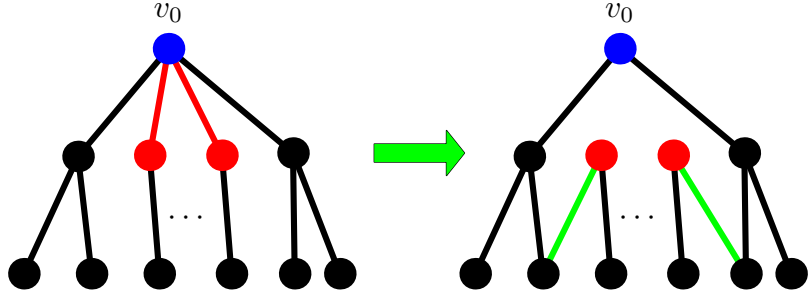}
\caption{The reconstruction of the tree with a vertex $v_0$ whose degree $d_{0} \geq 3$. For the original tree (i.e., the left-hand side of the figure), removing the links between vertex $v_0$ and its $d_{0}-2$ neighbors such that the degree of $v_0$ is $2$ and then connecting these neighbors to  $d_{0}-2$ leaf-vertices respectively, we obtain a new tree (i.e., the right-hand side of the figure) with less GMFRT.}
\label{fig:1}
\end{center}
\end{figure}

In fact, for any tree, let $k$ be the total number of leaf-vertices in the tree. If there is a vertex $v_0$ with degree $d_{0} \geq 3$ in this tree,  the degrees of other $n-k-1$ vertices are greater than or equal to $2$. Note that the sum of the degrees for all the vertices is $2(n-1)$. We have, $d_{0}+k+2(n-k-1)\leq 2(n-1).$ Therefore $k\geq d_{0}$.
Removing the links between vertex $v_0$ and its $d_{0}-2$ neighbors such that the degree of $v_0$ is $2$, and connectting the  $d_{0}-2$ old neighbors of $v_0$ to  $d_{0}-2$ vertices with  degree  $1$,  we obtain a new tree. The reconstruction of the tree is shown as Fig.~\ref{fig:1}. Let $N_1^{old}$, $N_2^{old}$ denote the number of vertices with degree $1$ and $2$ respectively in the old tree,   $N_1^{new}$, $N_2^{new}$ denote the number of vertices with degree $1$ and $2$ in the new tree respectively.  We have

 \begin{eqnarray}
N_1^{old}-N_1^{new}=d_{0}-2,\label{DEF1} \\
N_2^{old}-N_2^{new}=1-d_{0}\label{DEF2}.
\end{eqnarray}
Let $\overline{\mu}_{old}$, $\overline{\mu}_{new}$ denote the GMFRT for the two trees respectively.  We can obtain from Eqs.~(\ref{Def_GMFRT}), (\ref{EqPijC1}), ~(\ref{DEF1}) and ~(\ref{DEF2}) that~\cite{note1},
 \begin{eqnarray}
&&\overline{\mu}_{old}-\overline{\mu}_{new}\nonumber \\
&=&\frac{1}{n}[\frac{2(n-1)}{d_{0}}+ {2(n-1)}(d_{0}-2)-\frac{2(n-1)}{2}(d_{0}-1)]\nonumber \\
      &=& \frac{1}{n}[\frac{2(n-1)}{d_{0}}+ (n-1)(d_{0}-3)].\label{Def_Case1}
\end{eqnarray}
Thus  $\overline{\mu}_{old}-\overline{\mu}_{new}>0$ while $d_{0} \geq 3$. That is to say, if there is a vertex whose degree is more than $2$,  we can reconstruct the tree and obtain a tree with less GMFRT. Therefore the tree with the minimal GMFRT  has no vertex with degree  more than $2$.


\subsection{The trees with the maximal GMFRT}
\label{sec:3_2}
In this section, we find that the stars have the maximal GMFRT among all trees with the same numbers of vertices. For any star with $n$ $(n>2)$ vertices, all vertices has degree $1$ except for the central vertex whose degree is $n-1$. Thus we can obtain the maxmal GMFRT from Eqs. Eqs.~(\ref{Def_GMFRT}) and (\ref{EqPijC1})
  \begin{eqnarray}
\overline{\mu}_{max}&=&\frac{2(n-1)(n-1)+ 2}{n} \nonumber \\
      &=& 2n-4+\frac{4}{n}.
\end{eqnarray}
Now we will prove that the stars have the maximal GMFRT among all trees.

For any tree, if it is not a star, there are more than $2$ vertices with degree greater than or equal to $2$. Let $v_1$, $v_2$ be the two vertices with the highest degree of the tree and $d_{1}$, $d_{2}$ denote the  degree of the two vertices respectively. Without loss of generality, we assume that   $d_{1} \geq d_{2} \geq 2$.

\begin{figure}
\begin{center}
\includegraphics[scale=1.2]{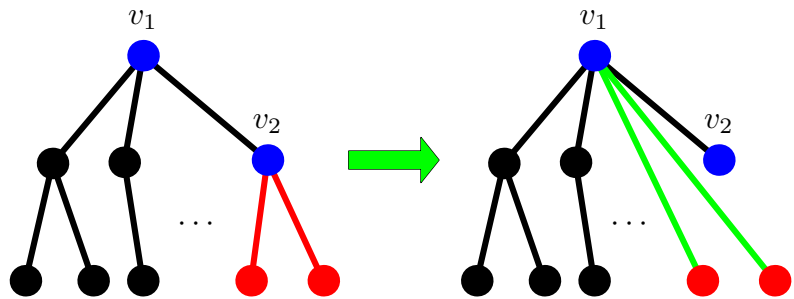}
\caption{The reconstruction of the tree with at least $2$ vertices (i.e., $v_1$,  $v_2$) whose degrees are greater than or equal to $2$. For the original tree (the left side of the figure), removing the links between vertex $v_2$ and its $d_{2}-1$ neighbors such that the degree of $v_2$ is $1$, connecting these neighbors to  vertex  $v_1$, we obtain a new tree (the right side of the figure) with larger GMFRT.}
\label{fig:2}
\end{center}
\end{figure}

As shown in Fig.~\ref{fig:2}, if we remove the links between vertex $v_2$ and its $d_{2}-1$ neighbors (except for the link between $v_2$ and $v_1$ if they are adjacent) and connect these neighbors of $v_2$ to  $v_1$, we obtain a new tree. In the new tree, the degree of $v_2$ is $1$ , the degree of $v_1$ is $d_{1}+d_{2}-1$,  the degrees of other vertices are the same as those of the old tree.

Let $\overline{\mu}_{old}$, $\overline{\mu}_{new}$  denote  the GMFRT for the two trees respectively.  We can obtain from Eqs.~(\ref{Def_GMFRT}) and (\ref{EqPijC1})  that,
 \begin{eqnarray}
\overline{\mu}_{old}-\overline{\mu}_{new}&=&\frac{ 2(n-1)}{n}[\frac{1}{d_{1}}+ \frac{1}{d_{2}}-\frac{1}{d_{1}+d_{2}-1} -1)] \nonumber \\
      &=& \frac{ 2(n-1)}{n}\cdot\frac{(d_{1}+d_{2})(1-d_{2})(d_{1}-1)}{d_{1}d_{2}(d_{1}+d_{2}-1)}.\nonumber
\end{eqnarray}
Thus  $\overline{\mu}_{old}-\overline{\mu}_{new}<0$ while $d_{1} \geq d_{2} \geq 2$. That is to say, if there are more than $2$ vertices with degree greater than or equal to $2$,  we can reconstruct a tree with larger GMFRT and the tree with the maximal GMFRT has only one vertex with degree  more than $2$, This kind of trees are just  stars.

 \section{The optimal  networks among all connected networks  with the same number of vertices and edges}
\label{sec:OP_N}

In this section, we  compare the  GMFRT  among all connected networks with $n$ vertices  $\{v_1, v_2,\cdots, v_n\}$  and $m$ edges, and then find the structures  who have the maximal GMFRT and the  structures who have the minimal GMFRT.

\subsection{The networks with the minimal GMFRT}
\label{sec:OP_N_1}
 In order to find the minimal  GMFRT among all connected networks with $n$ vertices , we should solve the following optimization equation
\begin{equation}
\label{LowerLimits}
 \left\{ \begin{array}{l}
  min ~~ \frac{1}{n} \sum_{i=1}^n\frac{1}{\pi_i} \\
  s.t.~~ \sum_{i=1}^n{\pi_i}=1,\pi_i>0(i=1,2,\ldots,n)
   \end{array} \right..
\end{equation}
Introducing a Lagrange multiplier $\lambda$ and studying the Lagrange function  defined by
\begin{equation}
F(\pi_1,\pi_2,\ldots,\pi_n,\lambda)= \frac{1}{n} \sum_{i=1}^n\frac{1}{\pi_i}+\lambda (\sum_{i=1}^n{\pi_i}-1).
\end{equation}
Setting the gradient $\nabla_{\pi_1,\pi_2,\ldots,\pi_n}F(\pi_1,\pi_2,\ldots,\pi_n,\lambda)=0$  yields the following equations:

\begin{equation}
\label{Lagrange_multi}
 \left\{ \begin{array}{l}
  \frac{\partial F}{\partial\pi_1}=-\frac{1}{n\pi_1^{2}}+\lambda=0  \\
  \frac{\partial F}{\partial\pi_2}=-\frac{1}{n\pi_2^{2}}+\lambda=0  \\
   ~~~~~~~\ldots \\
   \frac{\partial F}{\partial\pi_n}=-\frac{1}{n\pi_n^{2}}+\lambda=0
   \end{array} \right..
\end{equation}
which shows that the minimum of Eq.(\ref{LowerLimits}) is obtained while  $\pi_1=\pi_2=\ldots=\pi_n$. By using the constraint $\sum_{i=1}^n{\pi_i}=1$, we find
\begin{eqnarray}
    \pi_i=\frac{1}{n}, \quad  \quad i=1,2,\ldots,n.  
    \label{Distr_Lower}
\end{eqnarray}
 Substituting $\pi_i$ from Eq.(\ref{Distr_Lower}) in Eq.(\ref{LowerLimits}), we find the the lower bound of GMFRT is equal to $n$, which is the  number of vertices  for the graph.
 Note that $\pi_i=\frac{d_{i}}{\sum_{j=1}^nd_{j}} \quad (i=1,2,\ldots,n)$. Therefore, the minimal GMFRT is obtained while  $d_{1}=d_{2}=\cdots=d_{n}$. That is to say,  the networks, whose vertices  have the same degree, have the minimal  GMFRT among all connected networks with the same number of vertices. For example, both ring with $n$ vertices and $n-$clique (i.e. a complete graph with $n$ vertices and $\frac{n(n-1)}{2}$ edges) have the minimal  GMFRT among all networks with  $n$ vertices because all their vertices have the same degree.

 Note that we  want to find the networks with the minimal  GMFRT among all connected networks with $n$ vertices and $m$ edges. For some $m$, we can construct a network with all the vertices  have the same degree. Then it just has the minimal GMFRT among all connected networks with $n$ vertices and $m$ edges. However, there are also some $m$, for which we can not construct a network with all the vertices  have the same degree. For example, if  $m=n-1$, the networks are trees, we can  not construct a network with all the vertices  have the same degree. For this cases, we can only construct the closest one: almost all nodes have the same degree. Let $k=\frac{2m}{n}$ be the average of the degree for all the vertices. The degrees for vertices of the network with the minimal  GMFRT can only be either $\lfloor k\rfloor$  or $\lceil k \rceil$. Here $\lfloor k\rfloor$ is the largest integer which is less than or equal to $k$ and $\lceil k \rceil$ is the smallest integer which is greater than or equal to $k$. Further more, the number of vertices with degree $\lceil k \rceil$ would be $2m-\lfloor k\rfloor \times n$, and  the number of vertices with degree $\lfloor k\rfloor$  would be $n-2m+\lfloor k\rfloor \times n$.

\subsection{The networks with the maximal GMFRT}
\label{sec:OP_N_2}

In this section, we analyze and find  the networks  who have the maximal GMFRT among all connected graphs with $n$ vertices (i.e., $\{v_1, v_2,\cdots, v_n\}$ ) and $m$ edges. As found in Sec.~\ref{sec:OP_N_1}, the networks with all the vertices  have the same degree are the networks with the minimal GMFRT.
  On the contrary, the networks with the maximal GMFRT would be the graphs whose vertex degrees have the biggest difference. Here we present the method for  constructing the networks with the maximal GMFRT for different cases.  It is easy to know that $n-1\leq m\leq\frac{n(n-1)}{2}$~\cite{note2}.

I) For the case while $m=\frac{n(n-1)}{2}$, we can not construct any other graph except for the complete graph with $n$ nodes. The complete graph is just the graph with the maximal GMFRT. 

II) For the case while $\frac{n(n-1)}{2}>m>\frac{(n-1)(n-2)}{2}+1$, we can construct a complete graph with $n-1$ vertices $v_1, v_2,\cdots, v_{n-1}$ firstly, then connect the last vertex $v_n$ with $v_1, v_2,\cdots, v_{k}$, where $k=m-\frac{(n-1)(n-2)}{2}$ and $2\leq k\leq n-2$. We obtain a graph $G_o$ with $n$ vertices and $m$ edges. We can also find that
\begin{equation} \label{caseii}
d_{i}^{G_o}=
\left\{                 
  \begin{array}{ll}
   n-1 & i=1,2,\cdots, k \\
   n-2 & i=k+1,k+2,\cdots,n-1\\
   k   &  i=n
  \end{array}
\right.,                 
 \end{equation}
and $d_{min}^{G_o}=d_{n}^{G_o}=k$,
where $d_{i}^{G_o}$, $d_{min}^{G_o}$ denote the degree of vertex $v_i$, the minimal degree for  vertices of graph $G_o$ respectively.

Further more, we find, for any connected graph  $G$ with $n$ vertices and $m$ edges,
\begin{eqnarray}\label{ComGMFRT_G}
\overline{\mu}_{G_o}\geq \overline{\mu}_{G},
\end{eqnarray}
 where $\overline{\mu}_{G}$ is the GMFRT for graph $G$ and the proof of Eq.~(\ref{ComGMFRT_G}) is presented in Appendix A. Therefore $G_o$ is just the graph with the maximal GMFRT for the case  while $\frac{n(n-1)}{2}>m>\frac{(n-1)(n-2)}{2}+1$. 

III) For the case while $n-1<m\leq\frac{(n-1)(n-2)}{2}+1$, there is an integer $k$ $(1\leq k\leq n-3)$ such that $\frac{(n-k-1)(n-k-2)}{2}+k+1< m\leq \frac{(n-k)(n-k-1)}{2}+k$. The graph $G_o$ with the maximal GMFRT can be constructed in the following way. Firstly, we construct a complete graph with $n-k-1$ vertices $v_1, v_2,\cdots, v_{n-k-1}$,  then connect the vertex $v_{n-k}$ with $v_1, v_2,\cdots, v_{l}$, where $l=m-\frac{(n-k-1)(n-k-2)}{2}-k\geq 2$, finally, connect the vertex $v_{1}$ with $v_{n-k+1}, v_{n-k+2},\cdots, v_{n}$ respectively. We obtain a graph $G_o$ with $n$ vertices and $m$ edges. The degrees for vertices of graph $G_o$ satisfy
\begin{equation} \label{caseiii}
d_{i}^{G_o}=
\left\{                 
  \begin{array}{ll}
   n-1 & i=1 \\
   n-k-1 & i=2,3,\cdots, l \\
   n-k-2 & i=l+1,l+2,\cdots,n-k-1\\
   l   &  i=n-k\\
   1  & i=n-k+1,n-k+2,\cdots,n
  \end{array}
\right. .                
 \end{equation}

\begin{figure}
\begin{center}
\includegraphics[scale=1.5]{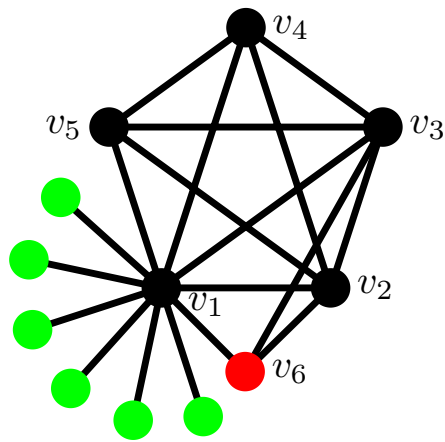}
\caption{The structure of the graph with the maximal GMFRT for the case while $n=12$, $m=19$.}
\label{fig:3}
\end{center}
\end{figure}

As proved  in Appendix B. for any connected graph  $G$ with $n$ vertices and $m$ edges,
\begin{eqnarray}\label{ComGMFRT_G_CASEIII}
\overline{\mu}_{G_o}\geq \overline{\mu}_{G}.
\end{eqnarray}
Therefore $G_o$ is just the graph with the maximal GMFRT for the case  while  $n-1<m\leq\frac{(n-1)(n-2)}{2}+1$. For example, for the case while $n=12$, $m=19$, we find that $\frac{(n-k-1)(n-k-2)}{2}+k+1< m\leq \frac{(n-k)(n-k-1)}{2}+k$, with $k=6$. Firstly, we construct a complete graph with $n-k-1=5$ vertices (i.e., $v_1, v_2,\cdots, v_{5}$),  then connect the vertex $v_{6}$ with $v_1, v_2, v_{3}$, finally, connect vertex $v_{3}$ with $v_{n-k+1}, v_{n-k+2},\cdots, v_{n}$ respectively. We obtain the graph with the maximal GMFRT, which is shown as  Fig.~\ref{fig:3}. 

IV) For the case while $m=n-1$, the graphs with $n$ vertices and $m$ edges are trees, the graph with the maximal GMFRT are stars. We have discussed this case in Sec.~\ref{sec:3_2}.

\section{Conclusion}
\label{sec:4}
We have found the  structures with the maximal GMFRT (i.e., the longest first return time) and the  structures with the minimal GMFRT (i.e., the shortest first return time) among all trees (and connected networks) with the same numbers of vertices. Our results show that, among all trees with the same number of vertices, trees with linear structure are the  structures with the minimal GMFRT and stars are the  structures with the maximal GMFRT.  We also find that, among all connected graphs with the same number of vertices, the graphs whose vertices  have the same degree, are the  structures with the minimal  GMFRT; and the graphs whose vertex degrees have the biggest difference, are the  structures with the minimal  GMFRT. We also present the methods for constructing the  graphs with the maximal GMFRT (or the minimal GMFRT), among all connected graphs with the same number of vertices and edges. If we use these graphs to mimic the state transition networks of extreme events, the  structures with the minimal GMFRT would be the optimal structures where the  recurrence time of the extreme events is the longest; the  structure with the minimal GMFRT would be the worst structures where the  recurrence time of the extreme events is the shortest. The finding of the optimal structures  would shed light on the control of extreme events. For example, we can lengthen the  recurrence time of the extreme events by increasing the difference of the vertex degrees for their state transition networks. The MFRT has close connection with the moments of first passage time, and the first passage time is useful indicator for transport efficiency of networks. The results obtained in this paper would be helpful for  controlling the transport efficiency of the networks.  In a opinion dynamic systems with random walking agents, we can control the fluctuation for the spreading of opinion by adjusting the network random walking agents move. If  the difference of the vertex degrees for the network is smaller,  the fluctuation for the spreading of opinion would be lower.

\ack{
This work is supported by the National Key R\&D Program of China (Grant No. 2018YFB0803604),  the National Natural Science Foundation of China (Grant No. 61873069, 61772147).

}

\appendix

\section{ Proof of  Eq.~(14) }
 \label{PR_E14}

 Note that, in graph $G_o$, vertices $v_1, v_2,\cdots, v_{n-1}$ have been connected with $\frac{(n-1)(n-2)}{2}$ edges and formed  a complete graph and $m=\frac{(n-1)(n-2)}{2}+d_{n}^{G_o}$.  We can not obtain any  graph $G$ with $n$ vertices, $m$ edges, and $d_{n}^{G}<d_{n}^{G_o}$. Therefore, for any graph  $G$ with $n$ vertices and $m$ edges,
  $$d_{min}^G \geq d_{min}^{G_o} =m-\frac{(n-1)(n-2)}{2}\geq 2,$$
  where  $d_{min}^{G}$ denotes the minimal degree for all nodes of graph $G$. 

For any connected graph $G$ with $n$ vertices and $m$ edges, if $d_{min}^G = d_{min}^{G_o}$,  $G$ and $G_o$ are isomorphic.They are  both the  the graph with the maximal GMFRT.
If $d_{min}^G > d_{min}^{G_o}$, we can reconstruct graph $G$  and obtain a graph  $G'$ which has  the same structure as graph $G_o$.
 Without loss of generality, let $v_n$ be the vertex whose degree is $d_{min}^G$ in graph $G$. Removing  links between vertex $v_n$ and its  $d_{min}^G -d_{min}^{G_o}$ neighbors, and then connect vertex $v_n$ with other vertices (e.g. $v_{i_1}, v_{i_2},\cdots, v_{i_l}$) whose degrees are less than $n-1$, we obtain a graph  $G'$. It easy to verify that $d_{min}^{G'} = d_{min}^{G_o}$ and the graph  $G'$ has the same structure as graph $G_o$. 
By comparing degree for nodes of the two graphs $G$ and $G'$, we find that
$$d_{min}^{G_o}= d_{n}^{G'}<d_{n}^G <n-1,$$
and $$d_{min}^{G_o}< d_{{i_j}}^G < d_{{i_j}}^{G'}=n-1, j=1,2,\cdots,l.$$
For any other vertex $v_i$ of graphs $G$ and $G'$,
$$d_{{i}}^G = d_{{i}}^{G'}.$$
Note that the change for degree of vertices  $v_{i_1}, v_{i_2},\cdots, v_{i_l}$ is just the change for vertex  $v_n$. We have
 $$d_{n}^G-d_{n}^{G'}=\sum_{j=1}^l \left( d_{{i_j}}^{G'}-d_{{i_j}}^G \right).$$
Therefore
 \begin{eqnarray}
\overline{\mu}_{G'}-\overline{\mu}_{G}&=&\frac{1}{n}\left[\frac{2m}{d_{n}^{G'}}-\frac{2m}{d_{n}^G}+\sum_{j=1}^l \left(\frac{2m}{ d_{{i_j}}^{G'}}-\frac{2m}{d_{{i_j}}^G}\right)\right]\nonumber \\
&=&\frac{2m}{n}\left[\frac{d_{n}^{G}-d_{n}^{G'}}{d_{n}^{G'}\times d_{n}^{G}}-\sum_{j=1}^l{\frac{d_{{i_j}}^{G'}-d_{{i_j}}^G} {d_{{i_j}}^{G'}\times d_{{i_j}}^G } }\right]\nonumber \\
&>&\frac{2m}{n}\left[\frac{d_{n}^{G}-d_{n}^{G'}}{d_{n}^{G'}\times d_{n}^{G}}-\sum_{j=1}^l{\frac{d_{{i_j}}^{G'}-d_{{i_j}}^G} {d_{n}^{G'}\times d_{n}^{G}} }\right]\nonumber \\
&=&0.
\end{eqnarray}
Thus Eq.~(\ref{ComGMFRT_G}) is obtained.

\section{ Proof of  Eq.~(16) }
 \label{PR_E16}

For the case while $n-1<m\leq\frac{(n-1)(n-2)}{2}+1$, there is an integer $k$ $(1\leq k\leq n-3)$ such that $\frac{(n-k-1)(n-k-2)}{2}+k+1< m\leq \frac{(n-k)(n-k-1)}{2}+k$. As discussed  in Sec.~\ref{sec:OP_N_2}, we can construct a graph $G_o$, which has $k$ vertices with degree $1$, a vertex with degree $l=m-\frac{(n-k-1)(n-k-2)}{2}-k\geq 2$, a vertex with degree $n-1$, and the degrees for other vertices are either $n-k-1$ or $n-k-2$. Here we will show that the graph $G_o$ has the  maximal GMFRT among all connected  graphs  with $n$ vertices and $m$ edges.

For any connected graph $G$  with $n$ vertices and $m$ edges, let $v_i$ $(i=1,2,\cdots,n)$ be its $n$ vertices and $d_{i}^G$ be the degree of vertex $v_i$. Without loss of generality, we assume that $d_{1}^G\geq d_{2}^G\geq \cdots \geq d_{n}^G$. If $d_{i}^G=d_{i}^{G_o}$ for any $(i=1,2,\cdots,n)$, graph $G$ is  isomorphic to $G_o$, they have the same GMFRT. If there is a vertex $v_{i_0}$ $(i_0=1,2,\cdots,n)$, such that $d_{{i_0}}^G\neq d_{{i_0}}^{G_o}$ and $d_{{i}}^G=d_{{i}}^{G_o}$ for $i>i_0$, we can adjust the graph $G$ and construct another graph $G'$ which has larger GMFRT than graph $G$. We will discuss the method for reconstructing the graph for different  $i_0$.

I) There is a vertex $v_{i_0}$ $(n-k+1 \leq i_0\leq n)$, such that $d_{{i_0}}^G\neq d_{{i_0}}^{G_o}$ and $d_{{i}}^G=d_{{i}}^{G_o}$ for $i>i_0$.

Note that $d_{i}^{G_o}=1$ for any $i\in [n-k+1, n]$. We can find that in graph $G$, the total number of vertices with degree $1$ is $n-i_0-1<k$. We remove the links between vertex $v_{i_0}$ and its $d_{{i_0}}^G-1$ neighbors such that the degree of $v_{i_0}$ is $1$, and for each of the $d_{{i_0}}^G-1$ old neighbors of $v_{i_0}$,we connect it to one of the vertices $\{v_1, v_2, \cdots, v_{i_0-1}\}$. We obtain a new connected graph $G'$ with $n$ vertices and $m$ edges.  We denote by  $v_{i_1}, v_{i_2},\cdots, v_{i_p}$ ($p\leq d_{{i_0}}^G-1$)  the vertices which had just received at least a link from the $d_{{i_0}}^G-1$ old neighbors of $v_{i_0}$. 

By comparing vertex degrees of the two graphs $G$ and $G'$, we find that
$$1=d_{{i_0}}^{G'}<d_{{i_0}}^G,$$
$$2\leq d_{{i_0}}^G\leq d_{{i_j}}^{G} < d_{{i_j}}^{G'}, j=1,2,\cdots,p,$$
and for any other vertex $v_i$ of graphs $G$ and $G'$,
$$d_{{i}}^G = d_{{i}}^{G'}.$$
Further more, $\sum_{j=1}^p\left(d_{{i_j}}^{G'}-d_{{i_j}}^{G}\right)=d_{{i_0}}^G-1$.
Therefore,
 \begin{eqnarray}
\overline{\mu}_{G'}-\overline{\mu}_{G}&=&\frac{1}{n}\left[{2m}-\frac{2m}{d_{{i_0}}^G}+\sum_{j=1}^p \left(\frac{2m}{ d_{{i_j}}^{G'}}-\frac{2m}{d_{{i_j}}^G}\right)\right]\nonumber \\
&=&\frac{2m}{n}\left[1-\frac{1}{d_{{i_0}}^G}-\sum_{j=1}^p{\frac{d_{{i_j}}^{G'}-d_{{i_j}}^G} {d_{{i_j}}^{G'}\times d_{{i_j}}^G } }\right]\nonumber \\
&>&\frac{2m}{n}\left[1-\frac{1}{d_{{i_0}}^G}-\sum_{j=1}^p{\frac{d_{{i_j}}^{G'}-d_{{i_j}}^G} {d_{{i_0}}^G\times (d_{{i_0}}^G+1) } }\right]\nonumber \\
&=&\frac{2m}{n}\left[1-\frac{1}{d_{{i_0}}^G}-{\frac{d_{{i_0}}^G-1} {d_{{i_0}}^G\times (d_{{i_0}}^G+1) } }\right]\nonumber \\
&>&\frac{2m}{n}\left[1-\frac{2}{d_{{i_0}}^G}\right]\nonumber \\
&\geq &0.
\end{eqnarray}

II) $i_0=n-k$, $d_{{i_0}}^G\neq d_{{i_0}}^{G_o}=l$ and $d_{{i}}^G=d_{{i}}^{G_o}=1$  for $i>n-k$.

We first show that $d_{{i_0}}^G>d_{{i_0}}^{G_o}=l$ for this case. Note that the graph $G$ has at least $k$ vertices (i.e., $v_i$ $(i=n-k+1, n-k+2,\cdots,n)$) with degree $1$. Removing the $k$ vertices (i.e., $v_i$ $(i=n-k+1, n-k+2,\cdots,n)$) from graph $G$, we obtain a sub-graph $G_1$ with $n-k$ vertices and $m-k$ edges. As $\frac{(n-k-1)(n-k-2)}{2}+k+1< m\leq \frac{(n-k)(n-k-1)}{2}+k$, we have
$$\frac{(n-k-1)(n-k-2)}{2}+1< m-k\leq \frac{(n-k)(n-k-1)}{2}.$$

Note that $\frac{(n-k-1)(n-k-2)}{2}$ is the total number of edges of a complete graph with $n-k-1$ vertices. The $m-k$ edges of $G_1$ are enough to construct a complete graph with $n-k-1$ vertices and then the minimal vertex degree for graph $G_1$ would be
$l=m-\frac{(n-k-1)(n-k-2)}{2}-k\geq 2$.  Therefore $d_{{i_0}}^G\geq d_{{i_0}}^{G_1}\geq l=d_{{i_0}}^{G_o}$. Recalling that $d_{{i_0}}^G\neq d_{{i_0}}^{G_o}$,
 we obtain $d_{{i_0}}^G> d_{{i_0}}^{G_o}=l$.

Now we adjust the graph $G$ and construct a new connected graph $G'$ with larger GMFRT than $G$. Removing the links between vertex $v_{i_0}$ and its $d_{{i_0}}^G-l$ neighbors such that the degree of $v_{i_0}$ is $l$, and for each of the $d_{{i_0}}^G-l$ old neighbors of $v_{i_0}$,we connect it to one of the vertices $\{v_1, v_2, \cdots, v_{i_0-1}\}$. We obtain a new connected graph $G'$ with $n$ vertices and $m$ edges.  We denote by  $v_{i_1}, v_{i_2},\cdots, v_{i_p}$ ($p\leq d_{{i_0}}^G-1$)  the vertices which had just received at least a link from the $d_{{i_0}}^G-l$ old neighbors of $v_{i_0}$. 

By comparing vertex degrees of the two graphs $G$ and $G'$, we find that
$$l=d_{{i_0}}^{G'}<d_{{i_0}}^G,$$
$$d_{{i_0}}^G\leq d_{{i_j}}^{G} < d_{{i_j}}^{G'}, j=1,2,\cdots,p,$$
and for any other vertex $v_i$ of graphs $G$ and $G'$,
$$d_{{i}}^G = d_{{i}}^{G'}.$$
Further more, $\sum_{j=1}^p\left(d_{{i_j}}^{G'}-d_{{i_j}}^{G}\right)=d_{{i_0}}^G-l$.
Therefore,
 \begin{eqnarray}
\overline{\mu}_{G'}-\overline{\mu}_{G}&=&\frac{1}{n}\left[\frac{2m}{ d_{{i_0}}^{G'}}-\frac{2m}{d_{{i_0}}^G}+\sum_{j=1}^p \left(\frac{2m}{ d_{{i_j}}^{G'}}-\frac{2m}{d_{{i_j}}^G}\right)\right]\nonumber \\
&=&\frac{2m}{n}\left[\frac{d_{{i_0}}^{G}-d_{{i_0}}^{G'}} {d_{{i_0}}^{G'}\times d_{{i_0}}^G }-\sum_{j=1}^p{\frac{d_{{i_j}}^{G'}-d_{{i_j}}^G} {d_{{i_j}}^{G'}\times d_{{i_j}}^G } }\right]\nonumber \\
&>&\frac{2m}{n}\left[\frac{d_{{i_0}}^{G}-l} {d_{{i_0}}^{G'}\times d_{{i_0}}^G }-\sum_{j=1}^p{\frac{d_{{i_j}}^{G'}-d_{{i_j}}^G} {d_{{i_0}}^{G'}\times d_{{i_0}}^G } }\right]\nonumber \\
&= &0.
\end{eqnarray}

III) There is a vertex $v_{i_0}$ $(1< i_0\leq n-k-1 )$, such that $d_{{i_0}}^G\neq d_{{i_0}}^{G_o}$ and $d_{{i}}^G=d_{{i}}^{G_o}$ for $i>i_0$. 

For this case,
\begin{equation} \label{caseiii}
d_{i}^{G}=
\left\{                 
  \begin{array}{ll}
   l   &  i=n-k\\
   1  & i=n-k+1,n-k+2,\cdots,n
  \end{array}
\right. .                
 \end{equation}
  We first show that there is no any edges between $v_{n-k}$ and $v_i$ $(i=n-k+1, n-k+2,\cdots,n)$ in graph $G$. On the contrary, if there are $p\geq 1$ edges between $v_{n-k}$ and $v_i$ $(i=n-k+1, n-k+2,\cdots,n)$ in graph $G$, we find that there are $\sum_{i=n-k}^n d_{{i}}^G-p =l+k-p$ edges connected with the $k+1$ vertices (i.e., $v_i$ $(i=n-k, n-k+1,\cdots,n)$) in graph $G$. If we remove the $k+1$ vertices (i.e., $v_i$ $(i=n-k, n-k+1,\cdots,n)$) and the $l+k-p$ edges from graph $G$, we obtain a sub-graph with $n-k-1$ vertices and $m-l-k+p=\frac{(n-k-1)(n-k-2)}{2}+p$ edges. The number of edges is more than the number of edges for a complete graph. It is impossible. Therefore there is no any edges between $v_{n-k}$ and $v_i$ $(i=n-k+1, n-k+2,\cdots,n)$ in graph $G$.

 Note that $d_{{i}}^G=d_{{i}}^{G_o}$ for $i>n-k$. If we remove the $k+1$ vertices (i.e., $v_i$ $(i=n-k, n-k+1,\cdots,n)$) from graph $G$ and $G_o$ at the same time. We obtain two sub-graphs with $n-k-1$ vertices and $\frac{(n-k-1)(n-k-2)}{2}$ from  graph $G$ and graph $G_o$. Both the two sub-graphs are complete graphs with $n-k-1$ vertices and have the same structure. Therefore the difference between $G$ and $G_o$ lies in how many vertices with degree $1$ are connected with vertices $v_1, v_2,\cdots, v_{n-k-1}$ respectively. For $G_o$, all vertices with degree $1$ are connected with $v_1$. But in graph $G$, there are also some vertices with degree $1$, which are connected with  vertices $v_{i}$ $(i=2,3,\cdots,i_0)$. By comparing vertex degree of the two graphs $G$ and $G_o$, we find that
  $$d_{{1}}^{G_o}>d_{{1}}^G, $$
$$ d_{{i}}^{{G_o}} < d_{{i}}^{G}, i=2,3,\cdots,i_0,$$
and for $i>i_0$,
$$d_{{i}}^{G_o} = d_{{i}}^{G}.$$
Further more, $\sum_{i=2}^{i_0}\left(d_{{i}}^{G}-d_{{i}}^{{G_o}}\right)=d_{{1}}^{G_o}-d_{{1}}^{G}$. Therefore,
 \begin{eqnarray}
\overline{\mu}_{G_o}-\overline{\mu}_{G}&=&\frac{1}{n}\left[\frac{2m}{ d_{{1}}^{G_o}}-\frac{2m}{d_{{1}}^G}+\sum_{i=2}^{i_0} \left(\frac{2m}{ d_{{i}}^{G_o}}-\frac{2m}{d_{{i}}^G}\right)\right]\nonumber \\
&=&\frac{2m}{n}\left[\sum_{i=2}^{i_0}{\frac{d_{{i}}^G -d_{{i}}^{G_o}} {d_{{i}}^{G_o}\times d_{{i}}^G } }-\frac{d_{{1}}^{G_o}-d_{{1}}^{G}} {d_{{1}}^{G_o}\times d_{{1}}^G }\right]\nonumber \\
&>&\frac{2m}{n}\left[\sum_{i=2}^{i_0}{\frac{d_{{i}}^G -d_{{i}}^{G_o}} {d_{{1}}^{G_o}\times d_{{1}}^G } }-\frac{d_{{1}}^{G_o}-d_{{1}}^{G}} {d_{{1}}^{G_o}\times d_{{1}}^G }\right]\nonumber \\
&= &0.
\end{eqnarray}

IV)  $d_{{1}}^G\neq d_{{1}}^{G_o}$ and $d_{{i}}^G=d_{{i}}^{G_o}$ for $n\geq i>1$. As two graphs $G$ and $G_o$ has the same number of vertices and edges, if $d_{{i}}^G=d_{{i}}^{G_o}$ for any $i$ ($n\geq i>1$), we have $d_{{1}}^G=d_{{1}}^{G_o}$. Therefore  it is impossible to find a graph $G$ which satisfies  $d_{{1}}^G\neq d_{{1}}^{G_o}$ and $d_{{i}}^G=d_{{i}}^{G_o}$ for any $i$ ($n\geq i>1$).


\begin{thebibliography}{100}
\bibitem{note1}
By comparing the vertex degrees of the old tree and new tree, there are $3$ differences: the degree for $v_0$, the number of vertexes with degree $1$ and the number of vertexes with degree $2$. The second line of Eq.~(\ref{Def_Case1}) shows the difference of MFRT between old tree and new tree. In the square brackets of the second line for Eq.~(\ref{Def_Case1}), the first term is the MFRT for walker startting from $v_0$ in old tree, the second term shows the difference of MFRT for vertices with degree $1$, the last term shows the  difference of MFRT for vertices with degree $2$. Because, in the new tree, the degree of $v_0$ is $2$, the MFRT for walker startting from $v_0$ in new tree is also recorded in the last term.
\bibitem{note2}
$m$ is less than or equal to the total number of edges for a complete graphs with $n$ vertices and $m$ is greater than or equal to the total number of edges for  a tree with $n$ vertices.
\end{thebibliography}
\end{document}